\documentclass{amsart}

\usepackage{amsmath,amsthm,amsfonts,amscd,amssymb}
\usepackage[pdfborder={0 0 0}]{hyperref}
\usepackage{enumerate}

\newcommand{\ZZ}{\mathbb{Z}}

\newtheorem{thm}{Theorem}[section]
\newtheorem{cor}[thm]{Corollary}
\newtheorem{lem}[thm]{Lemma}

\newtheorem{quest}{Question}

\theoremstyle{definition}

\theoremstyle{remark}
\newtheorem{rem}{Remark}[section]

\begin{document}

\title{On well-rounded ideal lattices - II}
\author[L. Fukshansky]{Lenny Fukshansky}
\author[G. Henshaw]{Glenn Henshaw}
\author[P. Liao]{Philip Liao} 
\author[M. Prince]{Matthew Prince}
\author[X. Sun]{Xun Sun}
\author[S. Whitehead]{Samuel Whitehead}\thanks{The authors were supported by a grant from the Fletcher Jones Foundation. The first author was also partially supported by a grant from the Simons Foundation (\#208969 to Lenny Fukshansky) and by the NSA Young Investigator Grant \#1210223.}

\address{Department of Mathematics, 850 Columbia Avenue, Claremont McKenna College, Claremont, CA 91711}
\email{lenny@cmc.edu}
\address{Department of Mathematics and Computer Science, Wesleyan University, Middletown, CT 06459}
\email{ghenshaw@wesleyan.edu}
\address{Department of Mathematics, Claremont McKenna College, Claremont, CA 91711}
\email{PLiao14@students.claremontmckenna.edu}
\address{Department of Mathematics, Harvey Mudd College, Claremont, CA 91711}
\email{mthwate@gmail.com}
\address{School of Mathematical Sciences, Claremont Graduate University, Claremont, CA 91711}
\email{foxfur\_32@hotmail.com}
\address{Department of Mathematics, Pomona College, Claremont, CA 91711}
\email{scw22009@mymail.pomona.edu}


\begin{abstract}
We study well-rounded lattices which come from ideals in quadratic number fields, generalizing some recent results of the first author with K. Petersen~\cite{lf:petersen}. In particular, we give a characterization of ideal well-rounded lattices in the plane and show that a positive proportion of real and imaginary quadratic number fields contains ideals giving rise to well-rounded lattices.
\end{abstract}

\maketitle

\def\A{{\mathcal A}}
\def\AA{{\mathfrak A}}
\def\B{{\mathcal B}}
\def\C{{\mathcal C}}
\def\D{{\mathcal D}}
\def\EE{{\mathfrak E}}
\def\F{{\mathcal F}}
\def\x{{\mathcal H}}
\def\I{{\mathcal I}}
\def\II{{\mathfrak I}}
\def\J{{\mathcal J}}
\def\K{{\mathcal K}}
\def\kk{{\mathfrak K}}
\def\L{{\mathcal L}}
\def\LL{{\mathfrak L}}
\def\M{{\mathcal M}}
\def\mm{{\mathfrak m}}
\def\MM{{\mathfrak M}}
\def\N{{\mathcal N}}
\def\O{{\mathcal O}}
\def\OO{{\mathfrak O}}
\def\PP{{\mathfrak P}}
\def\R{{\mathcal R}}
\def\PNR{{\mathcal P_N(\real)}}
\def\PMNR{{\mathcal P^M_N(\real)}}
\def\PdNR{{\mathcal P^d_N(\real)}}
\def\s{{\mathcal S}}
\def\V{{\mathcal V}}
\def\X{{\mathcal X}}
\def\Y{{\mathcal Y}}
\def\Z{{\mathcal Z}}
\def\H{{\mathcal H}}
\def\cee{{\mathbb C}}
\def\Nn{{\mathbb N}}
\def\pee{{\mathbb P}}
\def\que{{\mathbb Q}}
\def\QQ{{\mathbb Q}}
\def\real{{\mathbb R}}
\def\RR{{\mathbb R}}
\def\zed{{\mathbb Z}}
\def\ZZ{{\mathbb Z}}
\def\aaa{{\mathbb A}}
\def\ff{{\mathbb F}}
\def\HDelta{{\it \Delta}}
\def\kk{{\mathfrak K}}
\def\qbar{{\overline{\mathbb Q}}}
\def\kbar{{\overline{K}}}
\def\ybar{{\overline{Y}}}
\def\kkbar{{\overline{\mathfrak K}}}
\def\ubar{{\overline{U}}}
\def\eps{{\varepsilon}}
\def\ahat{{\hat \alpha}}
\def\bhat{{\hat \beta}}
\def\gt{{\tilde \gamma}}
\def\h{{\tfrac12}}
\def\be{{\boldsymbol e}}
\def\bei{{\boldsymbol e_i}}
\def\bc{{\boldsymbol c}}
\def\bm{{\boldsymbol m}}
\def\bk{{\boldsymbol k}}
\def\bi{{\boldsymbol i}}
\def\bl{{\boldsymbol l}}
\def\bq{{\boldsymbol q}}
\def\bu{{\boldsymbol u}}
\def\bt{{\boldsymbol t}}
\def\bs{{\boldsymbol s}}
\def\bv{{\boldsymbol v}}
\def\bw{{\boldsymbol w}}
\def\bx{{\boldsymbol x}}
\def\bX{{\boldsymbol X}}
\def\bz{{\boldsymbol z}}
\def\bwy{{\boldsymbol y}}
\def\bY{{\boldsymbol Y}}
\def\bL{{\boldsymbol L}}
\def\ba{{\boldsymbol a}}
\def\bb{{\boldsymbol b}}
\def\bet{{\boldsymbol\eta}}
\def\bxi{{\boldsymbol\xi}}
\def\bo{{\boldsymbol 0}}
\def\bone{{\boldsymbol 1}}
\def\bol{{\boldsymbol 1}_L}
\def\ep{\varepsilon}
\def\p{\boldsymbol\varphi}
\def\q{\boldsymbol\psi}
\def\rank{\operatorname{rank}}
\def\aut{\operatorname{Aut}}
\def\lcm{\operatorname{lcm}}
\def\sgn{\operatorname{sgn}}
\def\spn{\operatorname{span}}
\def\md{\operatorname{mod}}
\def\Norm{\operatorname{Norm}}
\def\dim{\operatorname{dim}}
\def\det{\operatorname{det}}
\def\Vol{\operatorname{Vol}}
\def\rk{\operatorname{rk}}
\def\ord{\operatorname{ord}}
\def\ker{\operatorname{ker}}
\def\div{\operatorname{div}}
\def\Gal{\operatorname{Gal}}
\def\GL{\operatorname{GL}}
\def\SNR{\operatorname{SNR}}
\def\WR{\operatorname{WR}}
\def\IWR{\operatorname{IWR}}
\def\scg{\operatorname{\left< \Gamma \right>}}
\def\swrh{\operatorname{Sim_{WR}(\Lambda_h)}}
\def\ch{\operatorname{C_h}}
\def\cht{\operatorname{C_h(\theta)}}
\def\scgt{\operatorname{\left< \Gamma_{\theta} \right>}}
\def\scgmn{\operatorname{\left< \Gamma_{m,n} \right>}}
\def\gat{\operatorname{\Omega_{\theta}}}
\def\mn{\operatorname{mn}}
\def\disc{\operatorname{disc}}

\section{Introduction and statement of results}
\label{intro}

Ideal lattices are important objects in number theory and discrete geometry, which have been extensively studied in a series of papers by Eva Bayer-Fluckiger and her co-authors in the 1990's and 2000's (see, for instance, \cite{bayer1}, \cite{bayer2}, \cite{bayer_nebe}). In this note, we consider the simplest kind of ideal lattices coming from quadratic number fields. Let $K$ be a quadratic number field, and let us write $\O_K$ for its ring of integers. Then $K=\que(\sqrt{D})$ (real quadratic) or $K=\que(\sqrt{-D})$ (imaginary quadratic), where $D$ is a positive squarefree integer. The embeddings $\sigma_1, \sigma_2 : K \to \cee$ can be used  to define the standard embedding $\sigma_K$ of $K$ into $\real^2$: if $K=\que(\sqrt{D})$, then $\sigma_K : K \to \real^2$ is given by $\sigma_K = (\sigma_1,\sigma_2)$; if $K=\que(\sqrt{-D})$, then $\sigma_2=\overline{\sigma_1}$, and $\sigma_K=(\Re(\sigma_1), \Im(\sigma_1))$, where $\Re$ and $\Im$ stand for real and imaginary parts, respectively. Each nonzero ideal $I \subseteq \O_K$ becomes a lattice of full rank in $\real^2$ under this embedding, which we will denote by $\Lambda_K(I) := \sigma_K(I)$. Such lattices are called planar {\it ideal lattices}.

Given a lattice  $\Lambda \subset \real^2$ of full rank with a basis $\ba_1,\ba_2$, we can write $A=(\ba_1 \ba_2)$ for the corresponding basis matrix, and then $\Lambda = A \zed^2$. The corresponding norm form is defined as 
$$Q_A(\bx) = \bx^t A^t A \bx,$$
and we say that the lattice is {\it integral} if the coefficient matrix $A^t A$ of this quadratic form has integer entries; it is easy to see that this definition does not depend on the choice of a basis. It is also easy to see that every ideal lattice is integral. We define $\det(\Lambda)$ to be $|\det(A)|$, again independent of the basis choice, and (squared) {\it minimum} or {\it minimal norm}
$$|\Lambda| = \min \{ \|\bx\|^2 : \bx \in \Lambda \setminus \{\bo\} \} = \min \{ Q_A(\bwy) : \bwy \in \zed^2 \setminus \{\bo\} \},$$
where $\|\ \|$ stands for the usual Euclidean norm. Then each $\bx \in \Lambda$ such that $\|\bx\|^2 = |\Lambda|$ is called a {\it minimal vector}, and the set of minimal vectors of $\Lambda$ is denoted by $S(\Lambda)$. A lattice $\Lambda$ is called {\it well-rounded} (abbreviated WR) if the set $S(\Lambda)$ contains two linearly independent vectors. These vectors form a basis for $\Lambda$, and we refer to such a basis as a {\it minimal basis}. WR lattices are important in discrete geometry and in a variety of optimization problems (see \cite{martinet}).

In this note, we study WR ideal lattices coming from quadratic number fields. The general investigation of WR ideal lattices has recently been started in~\cite{lf:petersen}, where, in particular, infinite families of WR ideal lattices coming from real and imaginary quadratic number fields have been constructed. We continue this investigation. Let us call an ideal $I$ in the ring of integers $\O_K$ of a quadratic number field $K=\que(\sqrt{\pm D})$ well-rounded (WR) if the corresponding planar lattice $\Lambda_K(I)$ is WR. We will say that a positive squarefree integer $D$ satisfies the {\it $\nu$-nearsquare condition} if it has a divisor $d$ with $\sqrt{\frac{D}{\nu}} \leq d < \sqrt{D}$, where $\nu > 1$ is a real number. Our main result is the following theorem.

\begin{thm} \label{ideal_IWR} If $D$ satisfies the 3-nearsquare condition, then the rings of integers of quadratic number fields $K=\que(\sqrt{\pm D})$ contain WR ideals; the statement becomes if and only if when $K=\que(\sqrt{-D})$. This in particular implies that a positive proportion (more than $1/5$) of real and imaginary quadratic number fields contain WR ideals, more specifically
\begin{equation}
\label{K_density}
\liminf_{N \to \infty} \frac{ \left| \left\{ K = \que(\sqrt{\pm D}) : K \text{ contains a WR ideal, } 0 < D \leq N \right\} \right|}{ \left| \left\{ K = \que(\sqrt{\pm D}) : 0 < D \leq N \right\} \right|} \geq \frac{\sqrt{3}-1}{2\sqrt{3}}.
\end{equation}
Moreover, for every $D$ satisfying the 3-nearsquare condition the corresponding imaginary quadratic number field $K=\que(\sqrt{-D})$ contains only finitely many WR ideals, up to similarity of the corresponding lattices, and this number is
\begin{equation}
\label{ideal_est_1}
\ll \min \left\{ 2^{\omega(D)-1}, \frac{2^{\omega(D)}}{\sqrt{\omega(D)}}  \right\},
\end{equation}
where $\omega(D)$ is the number of prime divisors of $D$ and the constant in the Vinogradov notation $\ll$ does not depend on $D$. 
\end{thm}

\begin{rem} \label{class_number} Two WR ideal lattices coming from  the same imaginary quadratic field are similar if and only if the corresponding ideals are in the same ideal class, which explains finiteness of the number of WR ideals in this case, as asserted in Theorem~\ref{ideal_IWR}. In fact, this number is strictly less than the class number $h_K$ unless $K=\que(\sqrt{-1})$ or $\que(\sqrt{-3})$, since $\O_K$ is WR if and only if $K$ is cyclotomic, as established in~\cite{lf:petersen}. We can then compare the bound of~\eqref{ideal_est_1} to estimates on the class number of imaginary quadratics. For instance, Siegel's original estimate~\cite{siegel_class_number} asserts that for $K=\que(\sqrt{-D})$,
\begin{equation}
\label{s-cl-n}
\log h_K \sim \log \sqrt{|\Delta|},
\end{equation}
where $\Delta$ is the corresponding fundamental discriminant:
$$\Delta = \left\{ \begin{array}{ll}
-D & \mbox{if $-D \equiv 1 (\md 4)$} \\
-4D & \mbox{if $-D \not\equiv 1 (\md 4)$,}
\end{array}
\right.$$
and so $h_K$ is about $O(\sqrt{D})$ as $D \to \infty$. On the other hand, the average order of $\omega(D)$ is $\log \log D$ (see~\cite{hardy}, \S\S22.11-22.13), and so the upper bound of~\eqref{ideal_est_1} is usually about $\frac{(\log D)^{\log 2}}{\sqrt{\log \log D}}$. Hence most ideal classes are not WR.
\end{rem}
\smallskip

The paper is structured as follows. In Section~\ref{iwr_section} we develop some notation, recalling a parameterization of integral well-rounded lattices in the plane, and prove a simple technical lemma which we use in our main argument. In Section~\ref{ideal} we demonstrate explicit constructions of WR ideals, in particular specifying infinite families of similarity classes of planar WR lattices containing ideal lattices, and then prove Theorem~1.1. In addition, we further discuss the case of real quadratic fields, as well as certain criteria for quadratic number fields to contain WR principal ideals at the end of Section~\ref{ideal}. We are now ready to proceed.
\bigskip

\section{Integral WR lattices in the plane}
\label{iwr_section}

We start with some notation, following~\cite{fletcher1}. An important equivalence relation on lattices is geometric similarity: two planar lattices $\Lambda_1, \Lambda_2 \subset \real^2$ are called {\it similar}, denoted $\Lambda_1 \sim \Lambda_2$, if there exists a positive real number $\alpha$ and a $2 \times 2$ real orthogonal matrix $U$ such that $\Lambda_2 = \alpha U \Lambda_1$. It is easy to see that similar lattices have the same algebraic  structure, i.e., for every sublattice $\Gamma_1$ of a fixed index in $\Lambda_1$ there is a sublattice $\Gamma_2$ of the same index in $\Lambda_2$ so that $\Gamma_1 \sim \Gamma_2$. If $\Lambda \subset \real^2$ is a full rank WR lattice, then its set of minimal vectors $S(\Lambda)$ contains 4 or 6 vectors, and this number is 6 if and only if $\Lambda$ is similar to the hexagonal lattice 
$$\H := \begin{pmatrix} 2 & 1 \\ 0 & \sqrt{3} \end{pmatrix} \zed^2$$
(see, for instance Lemma 2.1 of \cite{lf:petersen}). Any two linearly independent vectors $\bx,\bwy \in S(\Lambda)$ form a minimal basis. While this choice is not unique, it is always possible to select $\bx,\bwy$ so that the angle $\theta$ between these two vectors lies in the interval $[\pi/3,\pi/2]$, and any value of the angle in this interval is possible. From now on when we talk about a minimal basis for a WR lattice in the plane, we will always mean such a choice. Then the angle between minimal basis vectors is an invariant of the lattice, and we call it the {\it angle of the lattice} $\Lambda$, denoted $\theta(\Lambda)$; in other words, if $\bx,\bwy$ is any minimal basis for $\Lambda$ and $\theta$ is the angle between $\bx$ and $\bwy$, then $\theta = \theta(\Lambda)$ (see \cite{hex} for details and proofs of the basic properties of WR lattices in $\real^2$). In fact, it is easy to notice that two WR lattices $\Lambda_1,\Lambda_2 \subset \real^2$ are similar if and only if $\theta(\Lambda_1)=\theta(\Lambda_2)$ (see \cite{hex} for a proof).

The following parameterization of integral well-rounded (IWR) lattices is discussed in~\cite{fletcher1}. Let $\Lambda \subset \real^2$ be an IWR lattice, then
\begin{equation}
\label{cos_sin}
\cos \theta(\Lambda) = \frac{p}{q},\ \sin \theta(\Lambda) = \frac{r\sqrt{D}}{q}
\end{equation}
for some $r,q,D \in \zed_{>0}$, $p \in \zed_{\geq 0}$ such that
\begin{equation}
\label{prqD}
p^2+Dr^2=q^2,\ \gcd(p,q)=1,\ \frac{p}{q} \leq \frac{1}{2}, \text{ and } D \text{ squarefree},
\end{equation}
and so $\Lambda$ is similar to
\begin{equation}
\label{OprqD}
\Omega_D(p,q) := \begin{pmatrix} q & p \\ 0 & r\sqrt{D} \end{pmatrix} \zed^2.
\end{equation}
Moreover, for every $p,r,q,D$ satisfying \eqref{prqD}, $\Omega_D(p,q)$ is an IWR lattice with the angle $\theta(\Omega_D(p,q))$ satisfying \eqref{cos_sin}, and $\Omega_D(p,q) \sim \Omega_{D'}(p',q')$ if and only if $(p,r,q,D) = (p',r',q',D')$. In addition, if $\Lambda$ is any IWR lattice similar to $\Omega_D(p,q)$, then
\begin{equation}
\label{min_lattice}
\left| \Lambda \right| \geq \left| \frac{1}{\sqrt{q}} \Omega_D(p,q) \right|,
\end{equation}
where the lattice $\frac{1}{\sqrt{q}} \Omega_D(p,q)$ is also IWR. Due to this property, we call $\frac{1}{\sqrt{q}} \Omega_D(p,q)$ a {\it minimal} IWR lattice in its similarity class.

We say that an IWR planar lattice $\Lambda$ is of {\it type $D$} for a squarefree $D \in \zed_{>0}$ if it is similar to some $\Omega_D(p,q)$ as in \eqref{OprqD}. As discussed in~\cite{fletcher1}, the type is uniquely defined, i.e., $\Lambda$ cannot be of two different types. Moreover, a planar IWR lattice $\Lambda$ is of type $D$ for some squarefree $D \in \zed_{>0}$ if and only if all of its IWR finite index sublattices are also of type $D$. If this is the case, $\Lambda$ contains a sublattice similar to $\Omega_D(p,q)$ for every 4-tuple $(p,r,q,D)$ as in~\eqref{prqD}. Hence the set of planar IWR lattices is split into types which are indexed by positive squarefree integers with similarity classes inside of each type $D$ being in bijective correspondence with solutions to the ternary Diophantine equation $p^2+r^2D=q^2$. 
\smallskip

We also need the following counting functions, as introduced in~\cite{fletcher1}. For fixed positive integers $D$ and $r$, $D$ squarefree, define
\begin{equation}
\label{r_count}
f(r) = \left| \left\{ (p,q) \in \zed_{>0}^2 : q^2-p^2 = r^2D,\ \gcd(p,q) = 1,\ 0< \frac{p}{q} \leq \frac{1}{2} \right\} \right|,
\end{equation}
as well as
$$f_1(r) = \left| \left\{ (p,q) \in \zed_{>0}^2 : q^2-p^2 = r^2D,\ \gcd(p,q) = 1 \right\} \right|,$$
and
\begin{equation}
\label{f_21}
f_2(r) = \left| \left\{ (p,q) \in \zed_{>0}^2 : q^2-p^2 = r^2D,\ 0< \frac{p}{q} \leq \frac{1}{2} \right\} \right|.
\end{equation}
Notice that
\begin{equation}
\label{f_f1_f2}
f(r) \leq \min \{ f_1(r), f_2(r) \}.
\end{equation}
The function $f_1(r)$ is well-studied; in particular, Theorem 6.2.4 of \cite{mollin} implies that 
\begin{equation}
\label{f1_est}
f_1(r) \leq 2^{\omega(r^2D)-1} = 2^{\omega(rD)-1}.
\end{equation}
On the other hand, equation~(34) of~\cite{fletcher1} establishes that
\begin{equation}
\label{f_22}
f_2(r) = \left| \left\{ b \in \zed_{>0} : b \mid r^2D,\ r \sqrt{D} < b \leq r\sqrt{3D} \right\} \right|,
\end{equation}
and then inequality~(38) of~\cite{fletcher1} guarantees that
\begin{equation}
\label{f2_est}
f_2(r) \leq O \left( \frac{\tau(r^2D)}{\sqrt{\omega(r^2D)}} \right),
\end{equation}
where $\tau(r^2D)$ is the number of divisors of~$r^2D$. We now prove a lemma, which will be useful to us in Section~\ref{ideal}.

\begin{lem} \label{r_1} The equation $p^2+D=q^2$ with squarefree positive integer $D$ has an integral solution $(p,q)$ satisfying $p/q \leq 1/2$ if and only if 
\begin{equation}
\label{D_div}
D=d_1d_2 \text{ for some } d_1,d_2 \in \zed_{>0} \text{ with } d_1 < d_2,\ \sqrt{\frac{D}{3}} \leq d_1 < \sqrt{D}.
\end{equation}
If this is the case, then $\gcd(d_1,d_2)=1$, and if $(p,q)$ is such a solution, then $\gcd(p,q)=1$. Moreover, if $D$ is a positive even squarefree integer, then there are no solutions to $p^2+D=q^2$.
\end{lem}

\proof
Notice that $p^2+D=q^2$ has an integral solution if and only if it has a positive integral solution. The number of positive integral solutions $(p,q)$ with $p/q \leq 1/2$ is then given by $f_2(1)$ as in \eqref{f_21}, and so the equation has solutions if and only if $f_2(1) \geq 1$. Hence \eqref{f_22} with $r=1$ implies that this happens if and only if
$$D=d_1d_2 \text{ for some } d_1,d_2 \in \zed_{>0} \text{ with } d_1 < d_2,\ \sqrt{D} < d_2 \leq \sqrt{3D},$$
which is equivalent to \eqref{D_div}. Now, if \eqref{D_div} holds, then $\gcd(d_1,d_2)=1$ since $D$ is squarefree. Similarly, suppose $(p,q)$ is a solution with $\gcd(p,q)=g$, then $g^2 \mid D$, and so $g=1$.

Finally, suppose $2 \mid D = q^2-p^2$. Then $2 \mid q-p$ or $2 \mid q+p$, which means that 2 divides both, $q-p$ and $q+p$, and so $2^2 \mid D$, which contradicts $D$ being squarefree. Hence $p^2+D=q^2$ has no solutions.
\endproof
\bigskip

\section{Ideal WR lattices in the plane}
\label{ideal}

In this section we discuss ideal well-rounded lattices in the plane, proving Theorem~\ref{ideal_IWR}. We start by setting some standard notation, following \cite{lf:petersen} (see also \cite{buell} for a detailed exposition). As in Section~\ref{intro} above, let $K = \que(\sqrt{\pm D})$ be a real or imaginary quadratic number field, where $D$ is a positive squarefree integer, and let us write $\O_K$ for its ring of integers. We have $\O_K=\zed[\delta]$, where
\begin{equation}
\label{delta}
\delta = \left\{ \begin{array}{ll}
- \sqrt{D} & \mbox{if $K=\que(\sqrt{D})$, $D \not\equiv 1 (\md 4)$} \\
\frac{1-\sqrt{D}}{2} & \mbox{if $K=\que(\sqrt{D})$, $D \equiv 1 (\md 4)$} \\
- \sqrt{-D} & \mbox{if $K=\que(\sqrt{-D})$, $-D \not\equiv 1 (\md 4)$} \\
\frac{1-\sqrt{-D}}{2} & \mbox{if $K=\que(\sqrt{-D})$, $-D \equiv 1 (\md 4)$.}
\end{array}
\right.
\end{equation}
The embeddings $\sigma_1, \sigma_2 : K \to \cee$ are given by
$$\sigma_1(x+y\sqrt{ \pm D}) = x+y\sqrt{ \pm D},\ \sigma_2(x+y\sqrt{ \pm D}) = x-y\sqrt{ \pm D}$$
for each $x+y\sqrt{\pm D} \in K$, where $\pm$ is determined by whether $K$ is a real or an imaginary quadratic field, respectively. The number field norm on $K$ is defined by
$$\Nn(x+y\sqrt{\pm D}) = \sigma_1(x+y\sqrt{ \pm D}) \sigma_2(x+y\sqrt{ \pm D}) = \left( x+y\sqrt{ \pm D} \right) \left( x-y\sqrt{ \pm D} \right).$$
Now $I \subseteq \O_K$ is an ideal if and only if 
\begin{equation}
\label{I_abg}
I = \{ ax + (b+g\delta)y : x,y \in \zed \},
\end{equation}
for some $a,b,g \in \zed_{\geq 0}$ such that
\begin{equation}
\label{abg}
b < a,\ g \mid a,b,\text{ and } ag \mid \Nn(b+g\delta).
\end{equation}
Such integral basis $a,b+g\delta$ is unique for each ideal $I$ and is called the {\it canonical basis} for $I$. 

For each ideal $I$, we consider the corresponding planar ideal lattice $\Lambda_K(I) = \sigma_K(I)$ as defined in Section~\ref{intro}. An investigation of well-rounded ideal lattices has been initiated in \cite{lf:petersen} with some special attention devoted to the planar case. In particular, it has been established in \cite{lf:petersen} that the full ring of integers $\O_K$ is WR if and only if $K=\que(\sqrt{-1}), \que(\sqrt{-3})$. On the other hand, infinite families of real and imaginary quadratic fields with WR ideals have been constructed in \cite{lf:petersen}. Here we extend and generalize these constructions, showing that many similarity classes of planar IWR lattices contain ideal lattices.

Let $K=\que(\sqrt{\pm D})$ and $I \subseteq \O_K$ be an ideal with the canonical basis $a, b+g\delta$, as above. It is then easy to check that $\Lambda_K(I)$ has the following shape, which we record in a convenient form for the proof of our next lemmas:

\begin{trivlist}

\item If $K=\que({\sqrt{D}})$, $D \not\equiv 1 (\md 4)$, then
\begin{equation}
\label{B1}
\Lambda_K(I) = \begin{pmatrix} a & b-g\sqrt{D} \\ a & b+g\sqrt{D} \end{pmatrix} \zed^2 = \begin{pmatrix} a-b+g\sqrt{D} & b-g\sqrt{D} \\ a-b-g\sqrt{D} & b+g\sqrt{D} \end{pmatrix} \zed^2.
\end{equation}

\item If $K=\que({\sqrt{D}})$, $D \equiv 1 (\md 4)$, then
\begin{equation}
\label{B2}
\Lambda_K(I) = \begin{pmatrix} a & \frac{2b+g}{2} - \frac{g\sqrt{D}}{2} \\ a &\frac{2b+g}{2} + \frac{g\sqrt{D}}{2} \end{pmatrix} \zed^2 = \begin{pmatrix} \frac{2a-2b-g}{2} + \frac{g\sqrt{D}}{2}  & \frac{2b+g}{2} - \frac{g\sqrt{D}}{2} \\ \frac{2a-2b-g}{2} - \frac{g\sqrt{D}}{2} & \frac{2b+g}{2} + \frac{g\sqrt{D}}{2} \end{pmatrix} \zed^2.
\end{equation}

\item If $K=\que({\sqrt{-D}})$, $-D \not\equiv 1 (\md 4)$, then
\begin{equation}
\label{B3}
\Lambda_K(I) = \begin{pmatrix} a & b \\ 0 & -g\sqrt{D} \end{pmatrix} \zed^2 = \begin{pmatrix} a-b & b \\ g\sqrt{D} & -g\sqrt{D} \end{pmatrix} \zed^2.
\end{equation}

\item If $K=\que({\sqrt{-D}})$, $-D \equiv 1 (\md 4)$, then
\begin{equation}
\label{B4}
\Lambda_K(I) = \begin{pmatrix} a & \frac{2b+g}{2} \\ 0 & - \frac{g\sqrt{D}}{2} \end{pmatrix} \zed^2 = \begin{pmatrix} \frac{2a-2b-g}{2} & \frac{2b+g}{2} \\ \frac{g\sqrt{D}}{2} & - \frac{g\sqrt{D}}{2} \end{pmatrix} \zed^2.
\end{equation}

\end{trivlist}

\begin{lem}\label{lat_to_ideal} Let $D$ be a positive odd squarefree integer satisfying \eqref{D_div} of Lemma~\ref{r_1} and let $(p,q)$ be a solution to the equation $p^2+D=q^2$ with $p/q \leq 1/2$. Then the similarity class of $\Omega_D(p,q)$ contains ideal lattices. More specifically, let
\begin{equation}
\label{ab_D}
(a,b) = \left\{ \begin{array}{ll}
\left( p+q, \frac{p+q-1}{2} \right) & \mbox{if $D \equiv 1 (\md 4)$} \\
\left( 2(p+q), p+q \right) & \mbox{if $D \equiv 3 (\md 4)$},
\end{array}
\right.
\end{equation}
and define
\begin{eqnarray}
\label{IJ_ab}
& & I = I(p,q) := \{ ax + (b+\delta)y : x,y \in \zed \} \subset \O_{\que(\sqrt{D})} \nonumber \\
& & J = J(p,q) := \{ ax + (b+\delta)y : x,y \in \zed \} \subset \O_{\que(\sqrt{-D})}
\end{eqnarray}
for this choice of $(a,b)$, where $\delta$ is defined accordingly as in \eqref{delta}. Then $I,J$ are WR ideals in their respective quadratic number rings, and the ideal lattices $\Lambda_{\que(\sqrt{D})}(I)$, $\Lambda_{\que(\sqrt{-D})}(J)$ belong to the similarity class of $\Omega_D(p,q)$.
\end{lem}

\proof
It is easy to verify that $(a,b)$ as in \eqref{ab_D} and $g=1$ satisfy the conditions of \eqref{abg} with $\delta$ chosen respectively as in \eqref{delta}. Therefore $I$ and $J$ defined in \eqref{IJ_ab} are indeed ideals with the corresponding canonical basis $a,b+\delta$. Now a straight-forward calculation shows that the ideal lattices $\Lambda_{\que(\sqrt{D})}(I), \Lambda_{\que(\sqrt{-D})}(J)$ are WR with the corresponding minimal basis matrices being the second matrices in formulas \eqref{B1}-\eqref{B4}, respectively, and cosine of the angle of each such lattice being $p/q$. This completes the proof.
\endproof

In fact, Lemma~\ref{lat_to_ideal} allows for an additional observation on WR ideal lattices in the plane, which we record below.

\begin{cor} \label{all_in_type} Suppose that $\Gamma$ is an IWR lattice of type $D$, where $D$ is as in Lemma~\ref{lat_to_ideal}. Then $\Gamma$ contains IWR sublattices similar to ideal lattices coming from ideals in $\que(\sqrt{\pm D})$.
\end{cor}

\proof
As discussed in Section~\ref{iwr_section} above, $\Gamma$ contains sublattices similar to $\Omega_D(p,q)$ for any $p,r,q$ satisfying $p^2+r^2D=q^2$ with $\gcd(p,q)=1$ and $p/q \leq 1/2$. Since $D$ is as in the statement of Lemma~\ref{lat_to_ideal}, there must exist such $p,q$ with $r=1$. Then $\Gamma$ must contain IWR sublattices similar to $\Lambda_{\que(\sqrt{D})}(I)$ and to $\Lambda_{\que(\sqrt{-D})}(J)$ for $I,J$ as in \eqref{IJ_ab} for each such choice of $p,q$.
\endproof

Next we prove that the WR ideal lattices coming from imaginary quadratic fields constructed in Lemma~\ref{lat_to_ideal} are all that there are, up to similarity.

\begin{lem} \label{ideal_D} Let $D \in \zed_{>0}$ be squarefree and let $K=\que(\sqrt{- D})$ be such that there exists a WR ideal $I \subset \O_K$. Then $D$ must satisfy \eqref{D_div} of Lemma~\ref{r_1} and $\Lambda_K(I) \sim \Omega_D(p,q)$ for some $p,q$ so that $p^2+D=q^2$, $\gcd(p,q)=1$, $p/q \leq 1/2$.
\end{lem}

\proof
We start with a general remark that applies to any planar lattice $\Lambda$. Given a basis matrix $A$ for $\Lambda$, there must exist a change of basis matrix
\begin{equation}
\label{U_b_c}
U = \begin{pmatrix} s_1 & s_2 \\ s_3 & s_4 \end{pmatrix} \in \GL_2(\zed)
\end{equation}
so that $B=AU$ is a basis matrix for $\Lambda$ corresponding to a Minkowski-reduced basis (in case $\Lambda$ is WR, this is our minimal basis). 

We are now ready to start the proof with notation as in the statement of the lemma. Assume that the canonical basis for the ideal $I$ is $a,b+g\delta$, then $I=gI'$, where $I'$ has canonical basis $\frac{a}{g}, \frac{b}{g}+\delta$ and $\Lambda_K(I) \sim \Lambda_K(I')$. Hence we can assume without loss of generality that $g=1$.

First suppose that $-D \not\equiv 1 (\md 4)$, $K=\que(\sqrt{-D})$, then $\Lambda_K(I)$ is as in \eqref{B3} with $g=1$. Let $U$ as in \eqref{U_b_c} be the change of basis matrix from the first basis matrix in \eqref{B3} to a minimal basis matrix. Then $\Lambda_K(I)$ is a WR lattice with minimal basis matrix
\begin{equation}
\label{b_im_n1}
B = \begin{pmatrix} as_1+bs_3 & as_2+bs_4 \\ -s_3\sqrt{D} & -s_4\sqrt{D} \end{pmatrix}.
\end{equation}
Since
\begin{equation}
\label{sin_th}
\sin \theta(\Lambda_K(I)) = \frac{\det \Lambda_K(I)}{|\Lambda_K(I)|} = \frac{r\sqrt{D}}{q},
\end{equation}
where $\gcd(r,q)=1$, we immediately deduce from \eqref{B3} and \eqref{b_im_n1} that
$$r = \frac{a}{\gcd(a, (as_1+bs_3)^2+Ds_3^2)},$$
where
$$(as_1+bs_3)^2+Ds_3^2 = a(as_1^2+2bs_1s_3) + (b^2+D)s_3^2$$
is divisible by $a$, by \eqref{abg}, since $\Nn(b+g\delta) = b^2+D$ in this case. Therefore $r$ must be equal to 1, and so $\Lambda_K(I) \sim \Omega_D(p,q)$ for some $p,q$ so that $p^2+D=q^2$, $\gcd(p,q)=1$, $p/q \leq 1/2$.

Next suppose that $-D \equiv 1 (\md 4)$, $K=\que(\sqrt{-D})$, then $\Lambda_K(I)$ is as in \eqref{B4} with $g=1$. Let $U$ as in \eqref{U_b_c} be the change of basis matrix from the first basis matrix in \eqref{B4} to a minimal basis matrix. Then $\Lambda_K(I)$ is a WR lattice with minimal basis matrix
\begin{equation}
\label{b_im_1}
B = \begin{pmatrix} as_1+(b+1/2)s_3 & as_2+(b+1/2)s_4 \\ -s_3\sqrt{D}/2 & -s_4\sqrt{D}/2 \end{pmatrix}.
\end{equation}
Analogously to the argument above,
$$r = \frac{2a}{\gcd(2a, 4a^2s_1^2+4a(2b+1)s_1s_3+((2b+1)^2+D)s_3^2} = 1,$$
since $(2b+1)^2+D$ is divisible by $2a$, by \eqref{abg}, because $\Nn(b+g\delta) = \frac{1}{4} ((2b+1)^2+D)$ in this case. Hence again $\Lambda_K(I) \sim \Omega_D(p,q)$ for some $p,q$ so that $p^2+D=q^2$, $\gcd(p,q)=1$, $p/q \leq 1/2$. This completes the proof of the lemma.
\endproof

In the case of a real quadratic field the situation appears to be more complicated. We propose the following question.

\begin{quest} \label{D_real_quad} Do there exist real quadratic fields $K=\que(\sqrt{D})$ with positive squarefree $D$ not satisfying \eqref{D_div} of Lemma~\ref{r_1} so that $\O_K$ contains WR ideals?
\end{quest}

\noindent
Computational evidence suggests that the answer to this question is no, however at the moment we only have the following partial result in this direction.

\begin{lem} \label{ideal_D_real} Let $D \in \zed_{>0}$ be squarefree and let $K=\que(\sqrt{D})$ be such that there exists a WR ideal $I = \left< a,b+g\delta \right> \subset \O_K$, where $a,b+g \delta$ is the canonical basis for~$I$. Assume in addition that $a \mid 2D$, then $D$ must satisfy \eqref{D_div} of Lemma~\ref{r_1} and $\Lambda_K(I) \sim \Omega_D(p,q)$ for some $p,q$ so that $p^2+D=q^2$, $\gcd(p,q)=1$, $p/q \leq 1/2$.  In particular, if
\begin{enumerate}
\item $D \not\equiv 1 (\md 4)$ and $\min \{ a^2, b^2+D \} \geq 2ab$, or
\item $D \equiv 1 (\md 4)$ and $\min \left\{ a^2, \frac{1}{4} ((2b+1)^2+D) \right\} \geq 2a(b+1)$,
\end{enumerate}
then $a \mid 2D$.
\end{lem}

\proof
First notice, as in the proof of Lemma~\ref{ideal_D}, that $I=gI'$, where $I'$ has canonical basis $\frac{a}{g}, \frac{b}{g}+\delta$ and $\Lambda_K(I) \sim \Lambda_K(I')$. Hence we can assume without loss of generality that $g=1$.

It is easy to see that the condition $a \mid 2D$ is equivalent to $a \mid b^2+D$ if $D \not\equiv 1 (\md 4)$, and to $a \mid (2b+1)^2+D$ if $D \equiv 1 (\md 4)$. Indeed, if $D \not\equiv 1 (\md 4)$, then \eqref{abg} implies that
\begin{equation}
\label{nrm_cond_1}
a \mid \Nn(b+g\delta) = b^2-D,
\end{equation}
and so $a \mid 2D$ if and only if $a \mid b^2+D$; if $D \equiv 1 (\md 4)$, then \eqref{abg} implies that
\begin{equation}
\label{nrm_cond_2}
a \mid \Nn(b+g\delta) = \frac{1}{4} ((2b+1)^2-D),
\end{equation}
and so $a \mid 2D$ if and only if $a \mid (2b+1)^2+D$.

Now suppose that $D \not\equiv 1 (\md 4)$, $K=\que(\sqrt{D})$, then $\Lambda_K(I)$ is as in \eqref{B1} with $g=1$. Let $U$ as in \eqref{U_b_c} be the change of basis matrix from the first basis matrix in \eqref{B1} to a minimal basis matrix. Then $\Lambda_K(I)$ is a WR lattice with minimal basis matrix
\begin{equation}
\label{b_re_n1}
B = \begin{pmatrix} as_1+(b-\sqrt{D})s_3 & as_2+(b-\sqrt{D})s_4 \\ as_1+(b+\sqrt{D})s_3 & as_2+(b+\sqrt{D})s_4 \end{pmatrix}.
\end{equation}
Then we must have:
\begin{equation}
\label{vnn1}
a^2s_1^2+2abs_1s_3+(b^2+D)s_3^2=a^2s_2^2+2abs_2s_4+(b^2+D)s_4^2,
\end{equation}
and analogously to the arguments in the proof of Lemma~\ref{ideal_D} above in the imaginary case, we have
\begin{equation}
\label{Dr1}
r = \frac{a}{\gcd(a, a^2s_1^2+2abs_1s_3+(b^2+D)s_3^2)}.
\end{equation}
Hence if $a \mid 2D$, then $r=1$, and so $\Lambda_K(I) \sim \Omega_D(p,q)$ for some $p,q$ so that $p^2+D=q^2$, $\gcd(p,q)=1$, $p/q \leq 1/2$. We will now show that if (1) is satisfied, then $a \mid 2D$. Take $U$ to be the identity matrix. The positive definite binary quadratic norm form corresponding to the basis matrix $B$ in this case will be
$$Q_B(x,y) = 2a^2x^2+4abxy+2(b^2+D)y^2.$$
Since $\min \{ a^2, b^2+D \} \geq 2ab$, this form must be Minkowski reduced, which means that $B$ must be a minimal basis matrix. Since $\Lambda_K(I)$ is WR, a reduced basis for the lattice must consist of vectors of the same length, i.e we must have $a^2=b^2+D$, and hence $a \mid 2D$ by the argument above.

Next suppose that $D \equiv 1 (\md 4)$, $K=\que(\sqrt{D})$, then $\Lambda_K(I)$ is as in \eqref{B2} with $g=1$. Let $U$ as in \eqref{U_b_c} be the change of basis matrix from the first basis matrix in \eqref{B2} to a minimal basis matrix. Then $\Lambda_K(I)$ is a WR lattice with minimal basis matrix
\begin{equation}
\label{b_re_1}
B = \begin{pmatrix} as_1+ \frac{(2b+1-\sqrt{D})s_3}{2} & as_2+ \frac{(2b+1-\sqrt{D})s_4}{2} \\ as_1+ \frac{(2b+1+\sqrt{D})s_3}{2} & as_2+ \frac{(2b+1+\sqrt{D})s_4}{2} \end{pmatrix}.
\end{equation}
Then we must have:
\begin{eqnarray}
\label{vn1}
& & 4a^2s_1^2+8abs_1s_3+4as_1s_3+((2b+1)^2+D)s_3^2 \nonumber \\ 
& = & 4a^2s_2^2+8abs_2s_4+4as_2s_4+((2b+1)^2+D)s_4^2,
\end{eqnarray}
and analogously to the arguments in the proof of Lemma~\ref{ideal_D} above in the imaginary case, we have
$$r = \frac{a}{\gcd(a, 4a^2s_1^2+8abs_1s_3+4as_1s_3+((2b+1)^2+D)s_3^2)} = 1.$$
Hence if $a \mid 2D$, then $r=1$, and so again $\Lambda_K(I) \sim \Omega_D(p,q)$ for some $p,q$ so that $p^2+D=q^2$, $\gcd(p,q)=1$, $p/q \leq 1/2$. We will now show that if (2) is satisfied, then $a \mid 2D$. Again, take $U$ to be the identity matrix. The positive definite binary quadratic norm form corresponding to the basis matrix $B$ in this case will be
$$Q_B(x,y) = 2a^2x^2+4a(b+1)xy+ \frac{1}{2} ((2b+1)^2+D)y^2.$$
Since $\min \left\{ a^2, \frac{1}{4} ((2b+1)^2+D) \right\} \geq 2a(b+1)$, this form must be Minkowski reduced, which means that $B$ must be a minimal basis matrix. Since $\Lambda_K(I)$ is WR, a reduced basis for the lattice must consist of vectors of the same length, i.e we must have $4a^2=(2b+1)^2+D$, and hence $a \mid 2D$ by the argument above.
\endproof

In addition, we have the following finiteness result for the number of WR ideals in a fixed imaginary quadratic number field.

\begin{lem} \label{K_fin_num} Suppose that $D$ satisfies \eqref{D_div} of Lemma~\ref{r_1} and $K=\que(\sqrt{-D})$. Then $K$ contains only finitely many WR ideals, up to similarity of the corresponding lattices, and this number is
\begin{equation}
\label{ideal_est}
\ll \min \left\{ 2^{\omega(D)-1}, \frac{2^{\omega(D)}}{\sqrt{\omega(D)}}  \right\},
\end{equation}
where the constant in the Vinogradov notation $\ll$ does not depend on $D$.
\end{lem}

\proof
Lemma~\ref{r_1} guarantees that there exist integer pairs $(p,q)$ such that 
\begin{equation}
\label{pq1}
p^2+D=q^2 \text{ for some } p,q \text{ with } \gcd(p,q)=1, p/q \leq 1/2.
\end{equation}
Now Lemma~\ref{lat_to_ideal} guarantees that if $K=\que(\sqrt{-D})$, then there exists a WR ideal $I \subseteq \O_K$ with $\Lambda_K(I) \in \Omega_D(p,q)$ for each $p,q$ satisfying \eqref{pq1}, and Lemma~\ref{ideal_D} implies that {\it all} WR ideals in $\O_K$ correspond to solutions of \eqref{pq1}. Then the number of WR ideals in $K$, up to similarity of the lattices $\Lambda_K(I)$, is precisely the number of pairs $p,q$ as in \eqref{pq1}. This number is precisely $f(1)$ as defined in \eqref{r_count}, which is estimated by \eqref{f_f1_f2}. Now applying \eqref{f1_est} and \eqref{f2_est}, and noticing that for a squarefree integer $D$
$$\frac{\tau(D)}{\sqrt{\omega(D)}} = \frac{2^{\omega(D)}}{\sqrt{\omega(D)}},$$
we obtain \eqref{ideal_est}. 
\endproof

\proof[Proof of Theorem \ref{ideal_IWR}]
The first part of the theorem along with \eqref{ideal_est_1} follow from Lemmas~\ref{lat_to_ideal}, \ref{ideal_D}, and~\ref{K_fin_num} above, so it is only left to establish \eqref{K_density}. Define sets
\begin{equation}
\label{A_set}
\A = \left\{ D : D \in \zed_{>0} \text{ squarefree} \right\},
\end{equation}
and
\begin{equation}
\label{B_set}
\B_{\nu} = \left\{ D : D \in \zed_{>0} \text{ squarefree with a divisor } \frac{\sqrt{D}}{\nu} \leq d < \sqrt{D} \right\},
\end{equation}
where $\nu > 1$ is a real number. Define also
$$\A(N) = \left\{ D \in \A : D \leq N \right\},\  \B_{\nu}(N) = \left\{ D \in \B_{\nu} : D \leq N \right\},$$
for any $N \in \zed_{>0}$. To prove \eqref{K_density} we simply need to show that
\begin{equation}
\label{D_3}
\liminf_{N \to \infty} \frac{ \left| \B_{\sqrt{3}}(N) \right| }{ \left| \A(N) \right| } \geq \frac{\sqrt{3}-1}{2\sqrt{3}}.
\end{equation}
An analogue of \eqref{D_3} for integers that are not necessarily squarefree has been established in Theorem 4.4 of \cite{wr1}. We will now adapt the proof of Theorem 4.4 of \cite{wr1} to account for the squarefree condition.

Theorem 333 of \cite{hardy} implies that there exist absolute constants $c_1,c_2$ such that
\begin{equation}
\label{aN}
\frac{6N}{\pi^2} + c_1\sqrt{N} \leq |\A(N)| \leq \frac{6N}{\pi^2} + c_2\sqrt{N}.
\end{equation}
Now, following Section~4 of \cite{wr1}, we define
$$I_{\nu}(n) = \left\{n^2, n(n-1), \dots, n \left(n-\left[\left(\frac{\nu-1}{\nu}\right)n\right]\right)\right\},$$
for each $n \in \zed_{>0}$, and let
$$I'_{\nu}(n) = \left\{ m \in I_{\nu}(n) : m \text{ is squarefree} \right\}.$$
Suppose $n$ is prime, then
$$I'_{\nu}(n) = \left\{ nm : \left(n-\left[\left(\frac{\nu-1}{\nu}\right)n\right]\right) \leq m < n,\ m \text{ is squarefree} \right\},$$
and so, by \eqref{aN}
\begin{eqnarray}
\label{I_prime}
|I'_{\nu}(n)| & = & \left| \A(n) \right| - \left| \A \left(n-\left[\left(\frac{\nu-1}{\nu}\right)n\right]\right) \right| \nonumber \\
& \geq & \frac{6}{\pi^2} \left[\left(\frac{\nu-1}{\nu}\right)n\right] + (c_1-c_2) \sqrt{n}.
\end{eqnarray}
Since each $I'_{\nu}(n) \subseteq I_{\nu}(n)$ and $I_{\nu}(n) \cap I_{\nu}(m) = \emptyset$ when $\gcd(n,m)=1$, by part (ii) of Lemma~4.2 of \cite{wr1}, we conclude that $I'_{\nu}(n) \cap I'_{\nu}(m) = \emptyset$ when $\gcd(n,m)=1$. Notice also that
$$\bigcup_{n=1}^{[\sqrt{N}]} I'_{\nu}(n) \subseteq \B_{\nu}(N).$$
Now we adapt the argument in the proof of Lemma~4.3 of \cite{wr1}. Let $M=\pi(\sqrt{N})$, i.e. the number of primes up to $\sqrt{N}$. A result of Rosser and Schoenfeld (Corollary~1 on p.~69 of \cite{rosser}) implies that for all $\sqrt{N} \geq 17$,
\begin{equation}
\label{N_bound}
\frac{\sqrt{N}}{\log \sqrt{N}} < M < \frac{1.25506\ \sqrt{N}}{\log \sqrt{N}}
\end{equation}
Hence suppose that $N \geq 289$, and let $p_1, \dots, p_M$ be all the primes up to $\sqrt{N}$ in ascending order. Then $I_{\nu}(p_i) \cap I_{\nu}(p_j) = \emptyset$ for all $1 \leq i \neq j \leq N$.  Therefore, using \eqref{I_prime}, we obtain:
\begin{equation}
\label{low2}
|\B_{\nu}(N)| \geq \sum_{i=1}^M |I'_{\nu}(p_i)| \geq \frac{6(\nu-1)}{\pi^2 \nu} \sum_{i=1}^M p_i + (c_1-c_2) \sum_{i=1}^M \sqrt{p_i}.
\end{equation}
A result of R. Jakimczuk \cite{jakimczuk} implies that
\begin{equation}
\label{prime_bound}
\sum_{i=1}^M p_i > \frac{M^2}{2}\ \log^2 M.
\end{equation}
Notice also that
$$(c_1-c_2) \sum_{i=1}^M \sqrt{p_i} \geq - |c_1-c_2| M^{3/2},$$
and combining this observation with \eqref{N_bound}, \eqref{low2}, and \eqref{prime_bound}, we obtain:
\begin{equation}
\label{low}
|\B_{\nu}(N)| >  \frac{6(\nu-1)N}{2\pi^2 \nu} \left( 1 - \frac{\log \log \sqrt{N}}{\log \sqrt{N}} \right) ^2 - |c_1-c_2| N^{3/4} \left( \frac{1.25506}{\log \sqrt{N}} \right)^{3/2}.
\end{equation}
Notice also that as $N \to \infty$, \eqref{aN} implies that $|\A(N)| \leq \frac{(6+\eps)N}{\pi^2}$ for any $\eps > 0$, and combining this observation with \eqref{low}, we obtain:
\begin{equation}
\label{any_nu}
\liminf_{N \to \infty} \frac{ \left| \B_{\nu}(N) \right| }{ \left| \A(N) \right| } > \frac{6(\nu-1)}{2(6+\eps) \nu} \geq \frac{\nu-1}{2\nu},
\end{equation}
since the choice of $\eps$ is arbitrary. Now \eqref{D_3} follows by taking $\nu=\sqrt{3}$.
\endproof

Finally, we briefly discuss WR lattices coming from principal ideals. It is well known that every ideal in a ring of integers of a number field can be generated by at most two elements, and principal ideals play a very special role in algebraic number theory: they correspond to the identity element of the class group of a number field. It is therefore natural to ask whether principal ideals in a quadratic number ring can be WR? Corollary~2.4 of \cite{lf:petersen} implies that if $K=\que \left( \sqrt{-D} \right)$, then $\O_K$ contains principal WR ideals if and only if $D=1,3$. In the real quadratic case, the situation is again more complicated. We propose the following question.

\begin{quest} \label{D_principal} Do there exist real quadratic fields $K=\que(\sqrt{D})$ with positive squarefree $D \not\equiv 1 (\md 4)$ so that $\O_K$ contains principal WR ideals?
\end{quest}

\noindent
Computational evidence suggests that the answer to this question is no. On the other hand, there do exist $D \equiv 1 (\md 4)$ so that $\O_{\que(\sqrt{D})}$ contains a WR principal ideal $I$. In Table~\ref{table2} below we present a few examples of number fields $K=\que(\sqrt{D})$, $D \equiv 1 (\md 4)$, so that the class number $h_K=1$, which contain WR ideals. We present these ideals in terms of their canonical integral bases with $\delta$ as in \eqref{delta}.

\begin{center} 
\begin{table}[!ht]
\caption{Examples of WR ideals in $K=\que(\sqrt{D})$, $D \equiv 1 (\md 4)$, $h_K=1$} 
\begin{tabular}{|l|l|l|l|} \hline
$D$& WR ideals & Their similarity classes $p/q$ ($r=1$) \\
\hline
$21$ & $\left< 3, 1+\delta \right>, \left< 7, 3+\delta \right>$ & $2/5, 2/5$\\ \hline
$77$ & $\left< 7, 3+\delta \right>, \left< 11, 5+\delta \right>$ & $2/9, 2/9$\\ \hline
$133$ & $\left< 7, 3+\delta \right>, \left< 19, 9+\delta \right>$ & $6/13, 6/13$\\ \hline
$209$ & $\left< 11, 5+\delta \right>, \left< 19, 9+\delta \right>$ & $4/15, 4/15$\\ \hline
\end{tabular}
\label{table2}
\end{table}
\end{center}
\bigskip

{\bf Acknowledgment.} We would like to thank the Fletcher Jones Foundation-supported Claremont Colleges research experience program, under the auspices of which a large part of this work was done during the Summer of 2011. We are also grateful to the referee for a thorough reading and helpful comments; in particular, Remark~\ref{class_number} was suggested by the referee.
\bigskip

\bibliographystyle{plain}  
\bibliography{wr_ideal-2}    
\end{document}